\documentstyle[12pt]{article}
\if@twoside  m
\else
\fi
   \newcommand{\ie}{{\em i.e.}}
   \newcommand{\eg}{{\em e.g.}}
   \newcommand{\cf}{{\em cf. }}
   \newcommand{\etc}{{\em etc. }}
   \newcommand{\QED}{\mbox{\rule[-1.5pt]{6pt}{10pt}}}
   \newcommand{\rhs}{{\em rhs }}
   \newcommand{\restr}{\vert\hskip -5.5pt
            \phantom{\vert}^{\scriptscriptstyle \backslash}}
   \newcommand{\R}{I\!\!R}
   \newcommand{\CC}{{\cal C}}
   \newcommand{\EE}{{\cal E}}
   \newcommand{\KK}{{\cal K}}
   \newcommand{\OO}{{\cal O}}
   \newcommand{\WW}{{\cal W}}
   \newcommand{\eps}{\varepsilon}
   \newtheorem{claim}{Claim}[section]
   \newtheorem{theorem}[claim]{Theorem}
   \newtheorem{lemma}[claim]{Lemma}
   \newtheorem{remark}[claim]{Remark}
   \newtheorem{conjecture}[claim]{Conjecture}

\begin{document}
\vspace*{20mm}

\begin{flushleft}
{\Large\bf Bound--state asymptotic estimates for \\ window--coupled
Dirichlet strips and layers}
\vspace{8mm}

Pavel Exner$^{a,b}$ and Simeon A. Vugalter$^a$
\vspace{5mm}

a) Nuclear Physics Institute, Academy of Sciences,
25068 \v{R}e\v{z} near Prague, \\
b) Doppler Institute, Czech Technical University, B\v rehov\'a 7,
11519 Prague, \\
{\em exner@ujf.cas.cz, vugalter@ujf.cas.cz}
\end{flushleft}
\vspace{8mm}

\noindent
We consider the discrete spectrum of the Dirichlet Laplacian on a
manifold consisting of two adjacent parallel straight strips or
planar layers coupled by a finite number $\,N\,$ of windows in the
common boundary. If the windows are small enough, there is just one
isolated eigenvalue. We find upper and lower asymptotic bounds on the
gap between the eigenvalue and the essential spectrum in the planar
case, as well as for $\,N=1\,$ in three dimensions. Based on these
results, we formulate a conjecture on the weak--coupling asymptotic
behaviour of such bound states.
\vspace{8mm}

\section{Introduction}

There has been some interest recently to Laplacians on strips or
layers. Such a system is trivial when the manifold is straight and
the boundary conditions are translationally invariant, so there is a
natural separation of variables. On the other hand, the spectral
properties become nontrivial if the transverse modes are coupled,
which can be achieved, \eg, if the manifold is bent, locally
deformed, or coupled to another one \cite{ES,DE,BGRS,ESTV,EV1,EV2}.

The interest stems from two sources. On the physical side, such
operators with Dirichlet boundary conditions are used as models of
various mesoscopic semiconductor structures. The corresponding
solid--state literature is rather rich --- see \cite{DE,ESTV} for
some references. On the other hand, bound states in systems with open
geometries pose also mathematical questions such as the
weak--coupling limit, validity of the semiclassical approximation,
resonance scattering in such structures, \etc Some properties of them
can be seen numerically \cite{ESTV} while analytical proofs are
missing. Recall also that a closely related problem concerns Neumann
Laplacians, namely the existence of trapped modes in acoustic
waveguides \cite{ELV,DE}.

In a recent paper \cite{EV1} we have studied a pair of parallel
Dirichlet strips of widths $\,d_1, d_2\,$ coupled laterally through a
window of a width $\,2a\,$ in the common boundary; we have shown that
there are positive $\,c_1,\, c_2\,$ such that the gap between the
ground state and the threshold of the continuous spectrum can be
estimated as
   \begin{equation} \label{2-1}
-c_1 a^4 \,\le\, \epsilon(a) -\left(\pi\over d\right)^2 \le\,
-c_2 a^4
   \end{equation}
for any $\,a\,$ small enough. The numerical result of \cite{ESTV}
suggests that the true asymptotics is of the same type, but proving
this assertion and finding the coefficient in the leading term
remains an open problem.

The aim of the present paper is to generalize the above inequalities
to the case of a finite number of connecting windows and to a higher
dimension. In the following section we shall prove the bounds for a
pair of strips with $\,N\,$ windows. In Section~3 we formulate the
analogous problem for two layers and prove two--sided asymptotic
bounds for a single window shrinking to a point. Proofs rely in both
cases on variational estimates and follow the same basic strategy as
in \cite{EV1}. On the other hand, the existence of multiple windows
or the change in dimension require numerous modifications, which
prompts us to present the argument with enough details.

The upper and lower asymptotics bounds we are going to derive are in
each case of the same type differing just by values of the constants.
We are convinced that ground state has an asymptotic expansion and
its lowest--order is given by functions analogous to our bounds. This
conjecture is formulated in the concluding section. At the same time,
our present method does not allow to squeeze the bounds, or even to
come close to the true values as the Remark~\ref{optimality} below
illustrates.

\section{$N\,$ windows in dimension two}

Consider a straight planar strip $\,\Sigma:= \R\times [-d_2,d_1]\,$.
Given finite sequences $\,\CC\equiv\{x_k\}_{k=1}^N\,$ of mutually distinct
points of the $\,x$--axis and $\,A=\{a_k\}_{k=1}^N\,$ with
$\,a_k>0\,$, we denote $\,\WW_k:=[x_k\!-\!a_k,x_k\!+\!a_k]\,$ and set
$\,\WW:= \bigcup_{k=1}^n \WW_k\,$. Then we define $\,H(d_1,d_2;\WW)\,$
as the Laplacian on $\,L^2(\Sigma)\,$ subject to the Dirichlet condition at
$\,y=-d_2,d_1\,$ as well as at the $\,\R\setminus\WW\,$ part of the
$\,x$--axis; this operator coincides with the Dirichlet Laplacian
at the strip with the appropriate piecewise cut (see Fig.~1) defined
in the standard way \cite[Sec.XIII.15]{RS4}.
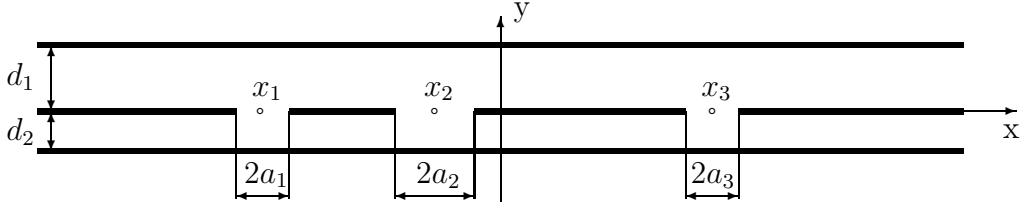
\begin{figure}
   \begin{picture}(120,80)
      \linethickness{2pt}
      \put(30,60){\line(1,0){350}}
      \put(30,35){\line(1,0){75}}
      \put(125,35){\line(1,0){40}}
      \put(195,35){\line(1,0){80}}
      \put(295,35){\line(1,0){85}}
      \put(30,20){\line(1,0){350}}
      \thinlines
      \put(205,1){\vector(0,1){70}}

      \put(108,8){$2a_1$}
      \put(105,3){\vector(1,0){20}}
      \put(125,3){\vector(-1,0){20}}
      \put(105,35){\line(0,-1){33}}
      \put(125,35){\line(0,-1){33}}
      \put(114,35){\circle{1.5}}
      \put(111,41){$x_1$}

      \put(173,8){$2a_2$}
      \put(165,3){\vector(1,0){30}}
      \put(195,3){\vector(-1,0){30}}
      \put(165,35){\line(0,-1){33}}
      \put(195,35){\line(0,-1){33}}
      \put(180,35){\circle{1.5}}
      \put(176,41){$x_2$}

      \put(277,8){$2a_3$}
      \put(275,3){\vector(1,0){20}}
      \put(295,3){\vector(-1,0){20}}
      \put(275,35){\line(0,-1){33}}
      \put(295,35){\line(0,-1){33}}
      \put(285,35){\circle{1.5}}
      \put(281,41){$x_3$}

      \put(380,35){\vector(1,0){20}}
      \put(35,35){\vector(0,1){25}}
      \put(35,60){\vector(0,-1){25}}
      \put(35,20){\vector(0,1){15}}
      \put(35,35){\vector(0,-1){15}}
      \put(18,45){$d_1$}
      \put(18,24){$d_2$}
      \put(210,71){y}
      \put(395,25){x}

   \end{picture}

\vspace{5mm}

\caption{Window--coupled planar waveguides}
   \end{figure}
Following the notation introduced in
\cite{EV1} we put $\,d:=\max\{d_1,d_2\}\,$ and $\,D:=d_1\!+d_2\,$. If
$\,d_1=d_2\,$, the operator decomposes into an orthogonal sum with
respect to the $\,y$--parity; the nontrivial part is unitarily
equivalent to the Laplacian on $\,L^2(\Sigma_+)\,$, where
$\,\Sigma_+:= \R\times [0,d]\,$, with the Neumann condition at window
part $\,\WW\,$ of the $\,x$--axis and Dirichlet at the remaining part
of the boundary; we denote it by $\,H(d;\WW)\,$. If the specification
is clear from the context, we will often denote the operator in
question simply as $\,H\,$.

We need a quantity to express the ``smallness" of the window set. We
define
   \begin{equation} \label{IW2}
I(\WW)\,:=\, \sum_{k=1}^N a_k|\WW_k|\,=\, 2\,\sum_{k=1}^N a_k^2\;;
   \end{equation}
then the result of \cite{EV1} generalizes to the present situation as
follows:

\begin{theorem} \label{2-N thm}
$\sigma_{\rm ess}(H(d_1,d_2;\WW))= [(\pi/d)^2\!,\infty)$. The discrete
spectrum is contained in $\,((\pi/D)^2\!,(\pi/d)^2)\,$, finite, and
non\-empty provided $\,\WW\ne\emptyset\,$. If $\,I(\WW)\,$ is
sufficiently small, $\sigma_{\rm disc}(H(d_1,d_2;\WW))\,$ consists of
just one simple eigenvalue $\,\epsilon(\WW) \le (\pi/d)^2$ and
there are positive $\,c_1,\, c_2\,$ such that
   \begin{equation} \label{2-N}
-c_1 I(\WW)^2 \,\le\, \epsilon(a) -\left(\pi\over d\right)^2 \le\, -c_2
I(\WW)^2
   \end{equation}
holds for any $\,I(\WW)\,$ small enough.
\end{theorem}

\noindent
{\em Proof: (a) The upper bound.} In the symmetric case,
$\,d_1=d_2\,$, the trial function will be chosen as $\,\psi=F+G\,$,
where
   \begin{equation} \label{trial F}
F(x,y)\,:=\, f_1(x)\chi_1(y)\,,
   \end{equation}
with
$$
f_1(x)\,:=\, \max\{\,\chi_{[x_1-a_1,x_N+a_N]}(x),\,
e^{-\kappa|x-x_1+a_1|},\, e^{-\kappa|x-x_N-a_N|}\}\,,
$$
and
   \begin{equation} \label{trial G}
G(x,y)\,:=\,\sum_{k=1}^N G_k(x,y)
   \end{equation}
with
   \begin{equation} \label{trial G_k}
G_k(x,y)\,:=\,{2\eta_k a_k\over |\WW|}\,\chi_{[x_k-a_k,x_k+a_k]}(x)\,
\cos\left(\pi (x\!-\!x_k)\over 2a_k\right)\, R_k(y)\,,
   \end{equation}
where $\,|\WW|:= 2\sum_{k=1}^N a_k\,$, and
   \begin{equation} \label{R}
R_k(y)\,:=\, \left\lbrace\; \begin{array}{lcc} e^{-\pi y/2a_k} & \quad
\dots\quad & y\in\left\lbrack 0,{d\over2}\,\right\rbrack \\ \\
2\left(1-{y\over d}\,\right) e^{-\pi d/4a_k} & \quad
\dots\quad & y\in\left\lbrack {d\over2},d\,\right\rbrack \end{array}
\right.
   \end{equation}
for $\,k=1,2,\dots,N\,$. As before $\,\chi_n(y)= \sqrt{2\over d}\,
\sin\left(\pi ny\over d\right)\,$, $\;n=1,2,\dots\,$, denote the
``transverse" eigenfunctions --- to be not confused with the
indicator function $\,\chi_M\,$ of a set $\,M\,$. Notice that as long
as we work with trial functions of $\,Q(H)\,$, the window smoothing
employed in \cite{EV1} is in fact not needed --- \cf\cite{RS4}.

The functional $\,L(\psi):= (H\psi, \psi)-\,\left(\pi\over
d\right)^2\|\psi\|^2$ can be expressed as
   \begin{equation} \label{trial L}
L(\psi)\,=\, \|\psi_x\|^2+ \|G_y\|^2 -\,
\left(\pi\over d\right)^2 \|G\|^2- 2\, {\pi\over d}\,
\sqrt{2\over d}\,\sum_{k=1}^N\, \int_{x_k-a_k}^{x_k+a_k} G_k(x,0)\,
dx\,.
   \end{equation}
Since $\,f_x,\, G_x\,$ have disjoint supports, we have
$\,\|\psi_x\|^2= \|F_x\|^2\!+ \sum_{k=1}^N \|G_{k,x}\|^2$, where
$\,G_{k,x}:= \partial_x G_k\,$. The $\,k$th term of the last sum
equals $\,\eta_k^2\pi^2 a_k |\WW|^{-2}\, \|R_k\|^2_{L^2(0,d)}\,$, and
$$
\|R_k\|^2_{L^2(0,d)}\,=\, {a_k\over\pi}\,+\, \left( {d\over 6}-
{a_k\over\pi} \right)\, e^{-\pi d/2a_k}\,<\, {a_k\over\pi}\,(1+\eps_1)
$$
for any $\,\eps_1>0\,$ and $\,a_k\,$ small enough. Obviously,
$\,\int_{x_k-a_k}^{x_k+a_k} G_k(x,0)\,dx=\, {8\over\pi}\, \eta_k
a_k^2|\WW|^{-1}\,$, and furthermore, a bound to $\,\|G_{k,y}\|^2$
follows from
$$
\|R'_k\|^2_{L^2(0,d)}\,=\, {\pi\over 4a_k}\,+\, \left( {2\over d}-
{\pi\over 4a_k} \right)\, e^{-\pi d/2a_k}\,<\, {\pi\over 4a_k}
$$
for $\,a_k<\pi d/8\,$, which means that $\,\|G_{k,y}\|^2<
\pi \sum_k \eta_k^2 a_k^2 |\WW|^{-2}\,$. Now we can put these
estimates together using $\,\|F_x\|^2= \kappa\,$; neglecting
the negative term $\,-\left(\pi\over d\right)^2 \|G\|^2\,$, we arrive
at the inequality
$$
L(\psi)\,<\, \kappa -\, {16\sqrt{2} \over
d^{3/2}}\, \sum_{k=1}^N {\eta_k a_k^2\over |\WW|}\,+\,
\pi(2+\eps_1)\sum_{k=1}^N {\eta_k^2 a_k^2\over |\WW|^2}\,.
$$
The sum of the last two terms at the \rhs is minimized by $\,-\,{2^7
\over \pi d^3(2+\eps_1)}\, \sum_k a_k^2\,$. To conclude the
argument, we have to estimate the trial function norm $\,\|\psi\|^2$
from below. The tail part is $\,\|\psi\|^2_{x\in\R\setminus\WW} =\,
\kappa^{-1}$, while the window contributes by
$$
\|\psi\|^2_{x\in\WW} \,\le \, 2 \|F\|^2_{x\in\WW}+
2\|G\|^2_{x\in\WW} \,=\, |x_N\!-\!x_1\!+\!a_N\!+\!a_1|\,
+ 4\, \sum_{k=1}^N {\eta_k^2 a_k^3\over |\WW|^2}\,
\|R_k\|^2_{L^2(0,d)}\,,
$$
so $\,\|\psi\|^2> (1\!-\!\eps_2) \kappa^{-1}\,$ holds for
any $\,\eps_2>0\,$ provided $\,|\WW|\,$ is small enough. Minimizing
the obtained estimate of $\,L(\psi)/\|\psi\|^2$ over $\,\kappa\,$, we
find
   \begin{equation} \label{trial bound}
{L(\psi)\over\|\psi\|^2}\,<\,-\, (1-\eps_2)^{-1} \left(\,{2^6 \over
\pi d^3(2+ \eps_1)}\, \sum_{k=1}^N a_k^2\, \right)^2
   \end{equation}
which yields the upper bound in (\ref{2-N}) for $\,d_1=d_2\,$. The
extension to the nonsymmetric case proceeds as for $\,N=1\,$; the
trial function is chosen in the above form for the wider channel,
while in the narrower one it is given by (\ref{trial G})
transversally rescaled.

\begin{remark} \label{optimality}
{\rm The bound can be improved, for instance, by replacing the
factorized form (\ref{trial G_k}) by a series, whose terms will be
products of the trigonometric basis in the window with the functions
$\,R_{k,n}(y)\,$ decaying as $\,\exp\{-{\pi ny\over 2a_k}\}\,$ around
$\,y=0\;$ (in the above estimate we used just the first term of such
a series). However, the gain is not large. To illustrate this fact,
take $\,N=1\,$ and $\,d=\pi\,$. The use of the series leads then to
the upper bound $\,\left(2a\over\pi\right)^4$ improving the
coefficient by $\,(\pi^2/8)^2\approx 1.52\,$. A comparison to the
numerically determined ground state \cite{ESTV} shows that the true
asymptotic behaviour should be $\,\approx (2.23\, a)^4$, so the
obtained $\,c_2\,$ is still two orders of magnitude off mark. The
reason is obviously that the wavefunction is affected by the window
outside the transverse ``window strip" as well. }
\end{remark}

Before proceeding to the lower bound, let us state some auxiliary
results:

\begin{lemma} \label{lemma 1}
Let $\,J[\phi]:= \int_a^b \left(\phi'(t)^2\!+ m^2\phi(t)^2 \right)\,
dt\,$ for $\,\phi\in C^2(a,b)\,$ with $\,\phi(a)=c_a\,$ (a fixed
number).  Given $\,m_0>0\,$, there is $\,\alpha_0>0\,$ such that
   \begin{equation} \label{interval bound}
J[\phi]\,\ge\, \alpha_0 m c_a^2
   \end{equation}
holds for all $\,m\ge m_0\,$.
\end{lemma}
{\em Proof:} The mimimum is obviously reached with $\,\phi'(b)=0\,$.
The corresponding Euler's equation is solved by $\,\phi_0(t)= d_1
e^{-mt}\!+d_2 e^{mt}\,$, where $\,d_1= c_a\left( e^{-ma}\!+
e^{m(a-2b)}\right)^{-1}\,$ and $\,d_2=d_1 e^{-2mb}\,$. Since
$\,m^{-1}c_a^{-2} \inf J(\phi)>0\,$ for any $\,m\ge m_0\,$, it is
sufficient to check that (\ref{interval bound}) remains valid as
$\,m\to\infty\,$; evaluating the functional for $\,\phi_0\,$ we find
$\,\lim_{m\to\infty} J(\phi)= mc_a^2\,$. \quad
\QED

\begin{lemma} \label{lemma 2}
Suppose that $\,\phi\,$ minimizes $\,J[\phi]:= \int_a^{2a}
\left(\phi'(t)^2\!+ p^2\phi(t)^2 \right)\, dt\,$ for positive
$\,a,\,p\,$ within $\,C^2(a,2a)\,$ with the boundary condition
$\,\phi(a)=c_a\,$; then
   \begin{equation} \label{exp bound}
|\phi(2a)| \,\le\, 2|c_a|\, e^{-pa}\,.
   \end{equation}
\end{lemma}
{\em Proof:} Assume for definiteness that $\,c_a>0\,$. By the
mentioned symmetry argument again, $\,\phi'(2a)=0\,$, and its explicit
form is $\,\phi(t)= c_a\,\cosh p(2a\!-\!t)/ \cosh pa\,$, which yields
$\,\phi(2a)\le 2c_a e^{-pa}$. \quad \QED
\vspace{3mm}

For the sake of completeness we reproduce also the following
assertion the proof of which is given in \cite{EV1}:

  \begin{lemma} \label{lemma 3}
Let $\,\phi\in C^2[0,d]\,$ with $\,\phi(0)=\beta\,$ and
$\,\phi(d)=0\,$. If $\,(\phi,\chi_1)=0\,$, then for every $\,m>0\,$
there is $\,d_0>0\,$ such that
   \begin{equation} \label{combined bound}
\int_0^d \phi'(t)^2 dt\,+\,\left(m\over a\right)^2 \int_0^a \phi(t)^2
dt\,-\, \left(\pi\over d\right)^2 \int_0^d \phi(t)^2 dt\,\ge\,
{d_0\beta^2 \over a}
   \end{equation}
holds for all $\,a\,$ small enough.
   \end{lemma}

{\em (b) Proof of Theorem~\ref{2-N thm}, continued:} The lower bound
is again the more difficult; however, we may restrict ourselves to
the symmetric case only because inserting an additional Neumann
boundary into the window we get a lower bound, and therefore we
consider in the following the spectrum of $\,H\equiv H(d;\WW)\,$.

We begin with a simple observation that it is sufficient to estimate
$\,L(\psi):= (H\psi, \psi)-\, \left(\pi\over d\right)^2\|\psi\|^2$
from below for all {\em real} $\,\psi\,$  of a core of $\,H$, say,
all $\,C^2$--smooth $\,\psi\in L^2(\Sigma_+)\,$ satisfying the
boundary conditions, since $\,H\,$ commutes with complex conjugation.
The main difficulty brought by the existence of multiple windows is
that we are no longer allowed to restrict ourselves to trial
functions symmetric with respect to the window centers. The strategy
we employ is to split from the beginning a part of the
kinetic--energy contribution to the functional, say, $\,{1\over 4}\,
\|\psi_x\|^2$, which will be at the end used to mend the problems
coming from the asymmetry, \ie, we begin with estimating
$\,L_0(\psi):= L(\psi)-\,{1\over 4}\,\|\psi_x\|^2$.

A trial function of the indicated set will be written in the form of
a Fourier series,
   \begin{equation} \label{Fourier}
\psi(x,y)\,=\, \sum_{n=1}^{\infty} c_n(x)\chi_n(y)
   \end{equation}
with smooth coefficients $\,c_n(x)=(\psi(x,\cdot),\chi_n)\,$, which
is uniformly convergent outside the windows, $\,x\not\in\WW\,$. We
split further the lowest transverse--mode coefficient by putting
   \begin{equation} \label{f_1}
f_1\,:=\, c_1\,-\, \sum_{k=1}^N \hat{f}_k\,,
   \end{equation}
where
   \begin{equation} \label{f_1 hat}
\hat f_k\,:=\, \left\lbrace \begin{array}{lll} c_k(x)\!-\!\alpha_k &
\quad \dots \quad & x\in[x_k\!-\!2a_k,x_k\!+\!a_k] \\ 0 & \quad
\dots \quad & {\rm otherwise} \end{array} \right.
   \end{equation}
with $\,\alpha_k:=c_1(x_k\!-\!2a_k)\,$, \ie, each one of the
functions $\,\hat f_k\,$ vanishes at the left endpoint of the
appropriate extended window; in contrast to \cite{EV1} we double the
left half of the window only. Writing the full trial function as
   \begin{equation}
\psi(x,y) \,=\, F(x,y)+G(x,y)\,, \qquad F(x,y)\,:=\,
f_1(x)\chi_1(y) \label{psi}\,,
   \end{equation}
we can cast the reduced energy functional into the form
   \begin{equation} \label{L}
L_0(\psi)\,=\, {3\over 4}\,\|\psi_x\|^2+ \|G_y\|^2 -\,
\left(\pi\over d\right)^2 \|G\|^2- \sum_{k=1}^N\, 2\alpha_k\,
{\pi\over d}\, \sqrt{2\over d}\, \int_{\WW_k} G(x,0)\,
dx\,.
   \end{equation}
Contributions to (\ref{L}) from different parts of the strip
$\,\Sigma_+\,$ will be estimated separately. The out--of--window part
consists of the sets
   \begin{eqnarray*}
\omega_1 \!&=&\! \{\, (x,y)\,:\; x\le x_1\!-\!a_1\,\}\,, \\
\omega_k \!&=&\! \{\, (x,y)\,:\; x_{k-1}\!+\!a_{k-1}\le x \le
x_k\!-\!a_k\,\}\,,\quad k=2,\dots,N\,, \\
\omega_{N+1} \!&=&\! \{\, (x,y)\,:\; x\ge x_N\!+\!a_N\,\}\,.
   \end{eqnarray*}
The expansion (\ref{Fourier}) yields
   \begin{eqnarray*}
\lefteqn{{1\over 4}\,\|\psi_x\|^2_{\omega_k}+ \|G_y\|^2_{\omega_k}
-\,\left(\pi\over d\right)^2 \|G\|^2_{\omega_k}} \\ \\ && =\,
{1\over 4}\,\sum_{n=1}^{\infty} \int_{\omega_k} c_n'(x)^2\, dx
\,+\,\sum_{n=1}^{\infty} \left(\pi\over d\right)^2
\left(n^2\!-1\right) \int_{\omega_k} c_n(x)^2\, dx \,,
   \end{eqnarray*}
and therefore
$$
{1\over 4}\,\|\psi_x\|^2_{\omega_k}+ \|G_y\|^2_{\omega_k}
-\,\left(\pi\over d\right)^2 \|G\|^2_{\omega_k} \,>\, \mu_0\,
\sum_{n=2}^{\infty} nc_n(x_k\!-\!a_k)^2
$$
with some $\,\mu_0>0\,$ follows from Lemma~\ref{lemma 1} (applied to
$\,c_n(-x)\,$) for $\,k=2,\dots,N\,$. The same inequality for
$\,k=1\,$ is derived as in \cite{EV1}; for the right tail we use just
the fact that the expression is positive so we can neglect it. Since
$\,\psi_x=G_x\,$ inside the (left extended) windows, we arrive at
the bound
   \begin{eqnarray} \label{L1}
L_0(\psi) &\!>\!& {1\over 2}\,\|\psi_x\|^2_{x\not\in\WW}+\,
\sum_{k=1}^N\, \bigg\lbrace\, {3\over 4}\,\|G_x\|^2_{x\in\WW_k} +
\|G_y\|^2_{x\in\WW_k} -\,\left(\pi\over d\right)^2
\|G\|^2_{x\in\WW_k} \nonumber \\ \\ &\!+\!& \mu_0\,
\sum_{n=2}^{\infty} nc_n(x_k\!-\!a_k)^2 -\, 2\alpha_k\, {\pi\over
d}\, \sqrt{2\over d}\, \int_{\WW_k} G(x,0)\, dx\, \bigg\rbrace\,.
\nonumber
   \end{eqnarray}
Our next goal is to estimate the contribution to $\,\|G_x\|^2$ from
the extended windows, $\,\EE_k:= [x_k\!-\!2a_k,x_k\!+\!a_k]\,$. In
distinction to the case $\,N=1\,$, however, even the lowest--mode
projection of $\,G\,$ may not vanish at the right endpoints of these
intervals, so the inequality (5.6) of \cite{EV1} has to be modified.
Fortunately, it is sufficient to change the coefficient: if a
function $\,\tilde G:\, \Sigma_+\to C^2(\Sigma_+)\,$ vanishes for
$\,x=x_k\!-\!2a_k\,$, the inequality (4.2) of \cite{EV1} in
combination with a symmetry argument imply
   \begin{equation} \label{zero end bound}
\|\tilde G_x\|^2_{x\in\EE_k}\,\ge\, \left(
\pi\over 6a_k\right)^2 \|\tilde G\|^2_{x\in\EE_k}\,.
   \end{equation}
To use this result we split the function by singling out the
projection of $\,G\,$ onto the first transverse mode,
   \begin{equation} \label{G}
G(x,y) \,=\, G_1(x,y)+G_2(x,y)\,, \qquad G_1(x,y)\,=\, \sum_{k=1}^N
\,\hat f_k(x)\chi_1(y)\,.
   \end{equation}
We have
$$
{1\over 2}\,\|\psi_x\|^2_{x\in\EE_k\setminus\WW_k}+\,{3\over 4}\,
\|G_x\|^2_{x\in\WW_k} \,\ge\, {1\over 2}\,\|G_x\|^2_{x\in\EE_k}  \,=\,
{1\over 2}\,\|G_{1,x}\|^2_{x\in\EE_k}+\, {1\over 2}\,
\|G_{2,x}\|^2_{x\in\EE_k}\,,
$$
and therefore
   \begin{eqnarray} \label{L2}
\lefteqn{L_0(\psi) \,>\, {1\over 2}\,\|\psi_x\|^2_{x\not\in\EE}+\,
\sum_{k=1}^N\, \bigg\lbrace\, {1\over 2}\,
\|G_{2,x}\|^2_{x\in\EE_k}\,+\, \|G_y\|^2_{x\in\WW_k}\,-\,
\left(\pi\over d\right)^2 \|G\|^2_{x\in\WW_k}} \nonumber \\ && \\
&& +\,  \mu_0\, \sum_{n=2}^{\infty} nc_n(x_k\!-\!a_k)^2 \,-\,
2\alpha_k\, {\pi\over d}\, \sqrt{2\over d}\, \int_{\WW_k}
G(x,0)\, dx\,+\,  {1\over 2}\, \left(\pi\over 6a_k\right)^2
\|G_1\|^2_{x\in\EE_k}  \bigg\rbrace\,. \nonumber
   \end{eqnarray}
with $\,\EE:= \bigcup_{k=1}^N\EE_k\,$. To proceed further we split
the function $\,G_2\,$ in the $\,k$--th extended window as
$\,G_2(x,y)= \hat G(x,y)+\Gamma(x,y)\,$, where
$$
\Gamma(x,y)\,:=\, \sum_{n=2}^{\infty} c_n(x_k\!-\!2a_k)\chi_n(y)\,.
$$
The second part is independent of $\,x\,$ while the first one
vanishes at left endpoint, so $\,G_{2,x}= \hat G_x\,$ may be
estimated by means of (\ref{zero end bound}) and the Schwarz
inequality as
   \begin{eqnarray} \label{G_2 bound}
\|G_{2,x}\|^2_{x\in\EE_k} &\!\ge\!&
\left(\pi\over 6a_k\right)^2 \|\hat G\|^2_{x\in\EE_k} \,\ge\,
\left(\pi\over 6a_k\right)^2 \|\hat G\|^2_{\Omega_k} \nonumber \\ \\
&\!\ge\!& {1\over 2}\left(\pi\over 6a_k\right)^2
\|G_2\|^2_{\Omega_k}\!- \left(\pi\over 6a_k\right)^2 \|\Gamma
\|^2_{\Omega_k}\,, \nonumber
   \end{eqnarray}
where we have denoted $\,\Omega_k:= \EE_k\times [0,a_k]\,$. To make
use of the last estimate we have to find  an upper bound to
$\,\|\Gamma\|^2_{\Omega_k}\,$. To this end we notice that
   \begin{description}
\item{\em (i)} instead of assuming $\,c_n\in C^2$, the lower bound
can be looked for in a wider class of $\,\psi\,$ with piecewise
continuous coefficients,
   \vspace{-1.8ex}
\item{\em (ii)} on the other hand, we may restrict ourselves to those
$\,\psi\,$ which satisfy for $\,x\in\EE_k\setminus \WW_k\,$ and
$\,n\ge 2\,$ the inequality
   \begin{equation} \label{subexponential decay}
|c_n(x)|\,\le\, c_n^{ex}(x)\,:=\, |c_n(a)|\, {\cosh\left( {\pi\over
d}\,\sqrt{n^2\!-1}\, (x\!-\!x_k\!+\!2a_k)\right) \over \cosh\left(
{\pi a_k\over d}\,\sqrt{n^2\!-1}\right)}\,.
   \end{equation}
To see that we split the trial function in analogy with \cite{EV1},
$$
\tilde\psi(x,y)\,:=\, \left\lbrace \begin{array}{lll}
\psi(x,y)-c_n(x)\chi_n(y) & \quad \dots \quad & x\in\EE_k
\setminus\WW_k \\ \psi(x,y) & \quad \dots \quad & {\rm otherwise}
\end{array} \right.
$$
The basic expression $\,L(\psi)/\|\psi\|^2$ can be then rewritten as
$$
\tilde L(\tilde\psi)-\left(\pi\over d\right)^2 \|\tilde\psi\|^2+
\sum_{k=1}^N \int_{\EE_k\setminus \WW_k} \left\lbrack\,
c_n'(x)^2\, dx + \left({\pi\over d}\,\sqrt{n^2\!-1}\,\right)^2
c_n(x)^2 \right\rbrack\, dx  \over
\|\tilde\psi\|^2 + \sum_{k=1}^N \int_{\EE_k\setminus \WW_k}
c_n(x)^2 dx \,,
$$
where $\,\tilde L(\tilde\psi):= \int_{\Sigma_+} \left(
|\tilde\psi_x|^2 +|\tilde\psi_y|^2 \right)(x,y)\, dx\,dy\,$.
We may assume only those $\,\psi\,$ for which the numerator is
negative; the part of its last term corresponding to the ``window
neighborhoods" is minimized by the hyperbolic function $\,c_n^{ex}\,$
of (\ref{subexponential decay}) (see the proof of Lemma~\ref{lemma
2}). It follows that replacing $\,c_n(x)^2\,$ by $\,\min\{c_n(x)^2,
c_n^{ex}(x)^2\}\,$ we can only get a larger negative number, while
the positive denominator can be only diminished.
   \end{description}
To estimate the norm of $\,\Gamma\,$ restricted to $\,\Omega_k\,$,
we adapt again the argument of \cite{EV1} and divide the series into
parts referring to small and large values of $\,y\,$, and employ,
respectively, the smallness of $\,\|\chi_n\restr [0,a]\|\,$ and the
bound (\ref{subexponential decay}). This yields
   \begin{eqnarray*}
&& \|\Gamma\|^2_{\Omega_k} \,=\, \int_{\EE_k} dx \int_0^{a_k} dy\,
\left( \sum_{n=2}^{\infty} c_n[2a_k]\chi_n(y) \right)^2 \\ \\
&\!\le\!& 6a_k \int_0^{a_k} \left( \sum_{n=2}^{[a_k^{-1}]+1}
c_n[2a_k]\chi_n(y) \right)^2\! dy \,+\, 6a_k \int_0^{a_k} \left(
\sum^{\infty}_{2\le n=[a_k^{-1}]+2} c_n[2a_k]\chi_n(y)
\right)^2 \!dy \\ \\
&\!\le\!& 24a_k \left( \sum_{n=2}^{[a_k^{-1}]+1} n^{-1}
c_n[a_k]^2 \int_0^{a_k} \chi_n(y)^2 dy \right) \left(
\sum_{n=2}^{[a_k^{-1}]+1} n \right) \\ \\
&\!+\!& 24a_k \left( \sum^{\infty}_{2\le n=[a_k^{-1}]+2} n
c_n[a_k]^2 \int_0^{a_k} \chi_n(y)^2 dy \right) \left(
\sum^{\infty}_{2\le n=[a_k^{-1}]+2} n^{-1} e^{-(2\pi
a_k/d)\sqrt{n^2-1}} \right)\,,
   \end{eqnarray*}
where $\,c_n[ja_k]:= c_n(x_k\!-\!ja_k)\,$ and $\,[\cdot]\,$ denotes
the entire part; in the the last step we have used the bound
$\,|c_n[2a_k]|<2|c_n[a_k]| \exp\left\lbrace -\,{\pi a_k\over d}
\sqrt{n^2\!-1} \right\rbrace\,$ which follows from Lemma~\ref{lemma
2}. In analogy with \cite{EV1}, this implies the existence of a
positive $\,C_k\,$ such that
   \begin{equation} \label{Gamma bound 2}
\|\Gamma\|^2_{\Omega_k}\,\le\, C_k a_k^2\, \sum_{n=2}^{\infty} n
c_n(x_k\!-\!a_k)^2\,.
   \end{equation}
From now on we consider continuous coefficient functions again. By
(\ref{G_2 bound}) we have
   \begin{eqnarray*}
\lefteqn{{1\over 2}\,\|G_{2,x}\|^2_{x\in\EE_k}+\,
\mu_0\,\sum_{n=2}^{\infty} n c_n(x_k\!-\!a_k)^2} \\ && \\ &&
 \ge\, \delta\left(\pi\over 12a_k\right)^2 \|G_2\|^2_{\Omega_k}\!-
{\delta\over 2}\, \left(\pi\over 6a_k\right)^2 \|\Gamma
\|^2_{\Omega_k}+\, \mu_0\,\sum_{n=2}^{\infty} n c_n(x_k\!-\!a_k)^2
   \end{eqnarray*}
for an arbitrary $\,\delta\in(0,1]\,$; if we choose the latter
sufficiently small, the sum of the last two terms is nonnegative for
each $\,k=1,\dots,N\,$ due to (\ref{Gamma bound 2}), so
   \begin{eqnarray} \label{L3}
L_0(\psi) &\!>\!& {1\over 2}\,\|\psi_x\|^2_{x\not\in\EE}+\,
\sum_{k=1}^N\, \bigg\lbrace\, \|G_y\|^2_{x\in\WW_k}
\,-\,\left(\pi\over d\right)^2 \|G\|^2_{x\in\WW_k}
\,+\, {m^2\over a_k^2}\, \|G_2\|^2_{\Omega_k} \nonumber \\ \\
&\!-\!& 2\alpha_k\, {\pi\over d}\, \sqrt{2\over d}\, \int_{\WW_k}
G(x,0)\, dx\,+\,  {1\over 2}\,\left(\pi\over 6a_k\right)^2
\|G_1\|^2_{x\in\EE_k} \bigg\rbrace\,, \nonumber
   \end{eqnarray}
where we have denoted $\,m:=\, {\pi\over 12}\, \sqrt{\delta}\,$.

Next we express the first term in the curly bracket using the
decomposition (\ref{G}), properties of the transverse base, and an
integration by parts,
$$
\|G_y\|^2_{x\in\WW_k}=\, \|G_{1,y}\|^2_{x\in\WW_k}+\,
\|G_{2,y}\|^2_{x\in\WW_k} -\,2\,{\pi\over d}\, \sqrt{2\over d}\,
\int_{\WW_k} \hat f_k(x)G(x,0)\, dx\,.
$$
As in \cite{EV1} we estimate the last term by the Schwarz inequality,
substitute into (\ref{L3}), neglect $\,\|G_{1,y}\|^2_{x\in\WW_k}$ as
well as
$$
{\pi^2\over 72 a_k^2}\, \|G_1\|^2_{x\in\EE_k}-\,
{\pi(\pi+\sqrt{2})\over d^2}\, \|G_{1,y}\|^2_{x\in\WW_k}
$$
which is positive for $\,a_k\,$ small enough, obtaining
   \begin{eqnarray} \label{L4}
L_0(\psi) &\!>\!& {1\over 2}\,\|\psi_x\|^2_{x\not\in\EE}+\,
\sum_{k=1}^N\, \bigg\lbrace\, \|G_{2,y}\|^2_{x\in\WW_k}
\,+\, {m^2\over a_k^2}\, \|G_2\|^2_{\Omega_k}\,-\,\left(\pi\over
d\right)^2 \|G\|^2_{x\in\WW_k}  \nonumber \\ \\
&\!-\!& {\pi\sqrt{2}\over d}\, \|G(\cdot,0)\|^2_{x\in\WW_k}
-\,2\alpha_k\, {\pi\over d}\, \sqrt{2\over d}\, \int_{\WW_k}
G(x,0)\, dx\, \bigg\rbrace\,. \nonumber
   \end{eqnarray}
By Lemma~\ref{lemma 3}, the sum of the first three terms in the curly
bracket is bounded from below by $\,{d_k\over a_k}\,
\|G(\cdot,0)\|^2_{x\in\WW_k}\,$ for some $\,d_k>0\,$. Since
$\,(d_k/2a_k) -(\pi\sqrt{2}/d)>0\,$ holds for $\,a_k\,$ small enough,
we have
$$
L_0(\psi) \,>\, {1\over 2}\,\|\psi_x\|^2_{x\not\in\EE}+\,
\sum_{k=1}^N \biggl\lbrace {d_k\over 2a_k}\,
\|G_2(\cdot,0)\|^2_{x\in\WW_k} \!-\,4\alpha_k\,{\pi\over d}\,
\sqrt{a_k\over d}\,\|G_2(\cdot,0)\|_{x\in\WW_k}\,
\biggr\rbrace\,,
$$
where we have employed again the Schwarz inequality. The $\,k$--th
term of the sum reaches its minimum w.r.t. the norm at $\,-\, {8\pi^2
\over d_kd^3}\, \alpha_k^2 a_k^2$. Returning to the original
functional and neglecting in the first term of the last estimate all
contributions except the one coming from the leftmost component of
$\,\R\setminus\EE\,$, we see that there is a positive $\,\gamma\,$
such that
   \begin{equation} \label{L5}
L(\psi) \,>\, {1\over 4}\,\|\psi_x\|^2+\, {1\over
2}\,\|\psi_x\|^2_{x<x_k\!-\! 2a_k} \,-\,\gamma\, \sum_{k=1}^N
\alpha_k^2 a_k^2
   \end{equation}
holds provided $\,|\WW|\,$ is small enough.

To conclude the proof, we denote $\,\ell_k:= x_k\!-\!2a_k\!-\!x_1
\!+\!2a_1\,$ and employ the identity
   \begin{equation} \label{alpha identity}
\sum_{k=1}^N \alpha_k^2 a_k^2\,=\, \sum_{k=1}^N \alpha_1^2 a_k^2
+\sum_{k=1}^N (\alpha_k^2\!-\alpha_1^2) a_k^2
   \end{equation}
together with the estimate
$$
{1\over 4}\,\|\psi_x\|^2\,\ge\, {1\over 4}\,\|c'_1\|^2 \,\ge\,
{1\over 4N}\, \sum_{k=1}^N\, {(\alpha_k\!-\alpha_1)^2\over \ell_k}\,.
$$
If $\, -\gamma(\alpha_k^2\!-\alpha_1^2) a_k^2+\, {1\over 4N\ell_k}
(\alpha_k\!-\alpha_1)^2 \ge 0\,$ holds for all $\,k=2\,\dots,N\,$,
the bound (\ref{L5}) reduces to
   \begin{equation} \label{L6}
L(\psi) \,>\, {1\over 2}\,\|\psi_x\|^2_{x<x_k\!-\! 2a_k}
\!-\,\gamma\alpha_1^2\,  \sum_{k=1}^N a_k^2\,.
   \end{equation}
On the other hand, suppose that the endpoint values satisfy
$\,\alpha_k\!-\!\alpha_1= \OO(a_k)\,$ as $\,a_k\to 0\,$ for
$\,k\in\KK\subset\{2,\dots,N\}\,$. In view of (\ref{alpha identity})
we have
$$
L(\psi) \,>\, {1\over 2}\,\|\psi_x\|^2_{x<x_k\!-\! 2a_k}
\!-\,\gamma\alpha_1^2\,  \sum_{k=1}^N a_k^2\,+\, \sum_{k\in\KK}\,
\left\lbrace\, {1\over 4N\ell_k}
(\alpha_k\!-\alpha_1)^2 -\,\gamma(\alpha_k^2\!-\alpha_1^2) a_k^2\,
\right\rbrace \;;
$$
however, the last term is $\,\OO\left(\sum_{k\in\KK} a_k^2
\right)\,$, so (\ref{L6}) is valid again with a smaller positive
coefficient in the last term. Since $\,\|\psi\|^2\ge
2\int_{-\infty}^{x_1\!-\!2a_1} c_1(x)^2 dx\,$, the quantity of
interest is bounded from below by
$$
{L(\psi)\over \|\psi\|^2} \,>\, {\int_{-\infty}^{x_1\!-\!2a_1}
c'_1(x)^2 dx \,-\,\gamma\alpha_1^2 I(\WW) \over
2\int_{-\infty}^{x_1\!-\!2a_1} c_1(x)^2 dx}\;.
$$
The \rhs is minimized by the function $\,c_1(x)= \alpha_1
e^{\kappa(x-x_1\!+2a_1)}$ which yields the value
$\,(\kappa^2/2)-\gamma I(\WW)\kappa\,$; taking the minimum over
$\,\kappa\,$ we find
   \begin{equation} \label{lower bound}
{L(\psi)\over \|\psi\|^2} \,>\, -\,\gamma^2 I(\WW)^2\,. \quad \QED
   \end{equation}

\section{Window--coupled layers}
\setcounter{equation}{0}

The setting of the three--dimensional problem is similar. We have a
straight layer, $\,\Sigma:= \R^2\times [-d_2,d_1]\,$, and a set
$\,\WW\subset \R^2\,$ which can be written as a finite union,
$\,\WW:= \cup_{k=1}^N \WW_k\,$, whose components are open, connected
sets of nonzero Lebesgue measure; without loss of generality we may
suppose they are mutually disjoint. Then we define $\,H(d_1,d_2;\WW)\,$
as the Laplacian on $\,L^2(\Sigma)\,$ obeying the Dirichlet condition at
the boundary of $\,\Sigma\,$, \ie, $\,y=-d_2,d_1\,$, as well as at
$\,\R^2\setminus\WW\,$. This operator coincides again with the
Dirichlet Laplacian \cite[Sec.XIII.15]{RS4} for the sliced layer the
two parts of which are connected through the window set $\,\WW\,$.
We use the same notation as above, $\,d:=\max\{d_1,d_2\}\,$ and
$\,D:=d_1\!+d_2\,$. The nontrivial part of the symmetric case,
$\,d_1=d_2\,$, reduces again to analysis of the Laplacian
$\,L^2(\Sigma_+)\,$, where $\,\Sigma_+:= \R^2\times [0,d]\,$, with
the Neumann condition at window part of the plane $\,y=0\,$ and
Dirichlet at the remaining part of the boundary; this operator will
be denoted as by $\,H(d;\WW)\,$.

Our main aim here is to prove a weak--coupling asymptotic estimate
for a pair of layers connected by a single window.

\begin{theorem} \label{3-1 thm}
$\sigma_{\rm ess}(H(d_1,d_2;\WW))= [(\pi/d)^2\!,\infty)$. The discrete
spectrum is contained in $\,((\pi/D)^2\!,(\pi/d)^2)\,$, finite, and
non\-empty provided $\,\WW\ne\emptyset\,$. Suppose further that
$\,N=1\,$ and $\,\WW=aM\,$ for an nonempty open set $\,M\,$ contained
in the unit ball $\,B_1\subset\R^2$. Then $\sigma_{\rm
disc}(H(d_1,d_2;aM))\,$ consists of just one simple eigenvalue
$\,\epsilon(aM) \le (\pi/d)^2$ for all $\,a\,$ small enough, and
there are positive $\,c_1,\, c_2\,$ such that
   \begin{equation} \label{3-1}
-\,\exp\left(-c_1 a^{-3}\right) \,\le\, \epsilon(a) -\left(\pi\over
d\right)^2 \le\, -\,\exp\left(-c_2 a^{-3}\right)\,.
   \end{equation}
\end{theorem}

\noindent
{\em Proof} is based again on variational estimates. {\em The upper
bound} in the symmetric case, $\,d_1=d_2\,$, employs the trial
function $\,\psi=F+\eta G\,$, where $\,F(x,y):= f_1(x)\chi_1(y)\,$
again with
   \begin{equation} \label{trial F_3}
f_1(x)\,:=\,\min\left\lbrace\, 1,\; {K_0(\kappa|x|)\over K_0(\kappa
a)} \right\rbrace\,,
   \end{equation}
and
   \begin{equation} \label{trial G_3}
G(x,y)\,:=\,\chi_{aM}(x)\phi_1(x)R(y)\,,
   \end{equation}
where $\,\phi_1^{(a)}\,$ is the ground--state eigenfunction,
$\,\|\phi_1^{(a)}\|=1\,$, of the operator $\,-\Delta_D^{aM}\,$
corresponding to the positive eigenvalue $\,\mu_1(a)=
\mu_1(1)a^{-2}\,$, and
   \begin{equation} \label{R_3}
R(y)\,:=\, \left\lbrace\; \begin{array}{lcc}
e^{-\sqrt{\mu_1(a)}\,y} & \quad \dots\quad & y\in\left\lbrack
0,{d\over2}\,\right\rbrack \\ \\
2\left(1-{y\over d}\,\right) \exp\left(-\,{d\over
2}\,\sqrt{\mu_1(a)}\right) & \quad \dots\quad & y\in\left\lbrack
{d\over2},d\,\right\rbrack \end{array} \right.
   \end{equation}
Using $\,-\chi''_1= (\pi/d)^2\chi_1\,$, a simple integration by
parts, and the fact that the vector functions $\,\nabla f_1\,$ and
$\,\nabla \phi_1^{(a)}\,$ have disjoint supports, we can express the
reduced energy functional $\,L(\psi):= (H\psi, \psi)-\,\left(\pi\over
d\right)^2\|\psi\|^2$ as
   \begin{eqnarray} \label{trial L_3}
L(\psi) \!&=&\! \|\nabla f_1\|^2_{L^2(\R^2)} +\eta^2 \left(
\mu_1(a) -\left(\pi\over d\right)^2 \right) \|R\|^2_{L^2(0,d)}
\nonumber \\ \!&-&\! \eta^2 \|R'\|^2_{L^2(0,d)} -2\eta \chi'_1(0)
\int_{aM} \phi_1^{(a)}(x)\, dx\,,
   \end{eqnarray}
where the negative term in the bracket can be, of course, neglected.
The second and the third term at the \rhs can be estimated in analogy
with \cite{EV1},
$$
\mu_1(a) \|R\|^2_{L^2(0,d)} -\eta^2 \|R'\|^2_{L^2(0,d)}\,<\,
{\sqrt{\mu_1(1)}\over 2a}\, (2\!+\eps_1)
$$
for a fixed $\,\eps_1>0\,$ and any $\,a\,$ small enough. In a similar
way, the last term equals $\,-2\eta\chi'_1(0)Ca\,$, where $\,C:=
\int_M \phi_1^{(1)}(x)\,dx\,$. Finally, the first one can be
evaluated by means of \cite[9.6.26]{AS},\cite[1.12.3.2]{PBM},
$$
K_0(\kappa a)^2 \|\nabla f_1\|^2_{L^2(\R^2)}\,=\, 2\pi\left\lbrack\,
{1\over 2}\,\kappa^2a^2 K'_1(\kappa a)^2 -\,{1\over 2}\left(\kappa^2
a^2\!+1\right) K_1(\kappa a)^2 \,\right\rbrack\,.
$$
Using $\,-K'_1(\xi)=K_0(\xi)+\xi^{-1}K_1(\xi)\,$ in combination with
the asymptotic expressions $\,K_0(\xi)=-\ln\xi+\OO(1)\,,\; K_1(\xi)=
\xi^{-1}\!+ \OO(\ln\xi)\,$, we find
$$
\|\nabla f_1\|^2_{L^2(\R^2)}\,<\, -\: {2\pi(1+\eps_2)\over
\ln\kappa a}
$$
for a fixed $\,\eps_2\,$ and $\,a\,$ small enough. Substituting these
estimates into (\ref{trial L_3}) and taking a minimum over $\,\eta\,$
we arrive at the bound
   \begin{equation} \label{upper L_3}
L(\psi)\,<\, -\: {2\pi(1+\eps_2)\over \ln\kappa a} \,-\,
{2\chi'_1(0)^2 C^2\over (2+\eps_1)\sqrt{\mu_1(1)}}\, a^3\,.
   \end{equation}
It remains to find a lower bound to
$$
\|\psi\|^2 \,\ge\, \|\psi\|^2_{|x|\ge a} -2\|F\|^2_{|x|\le a}
-2\eta^2 \|F\|^2_{|x|\le a} \,=\, \|\psi\|^2_{|x|\ge a} -2\pi a^2
-2\eta^2\|R\|^2_{L^2(0,d)}\,.
$$
The last term is $\,\OO(a)\,$, while the first one can be expressed
as
$$
K_0(\kappa a)^2 \|F\|^2_{|x|\ge a}\,=\, \pi a^2 \left\lbrack\,
K_1(\kappa a)^2\! -K_0(\kappa a)^2\,\right\rbrack \,=\, {\pi\over
\kappa^2}\,+\, \OO(a^2\ln\kappa a)\;;
$$
using the asymptotic behaviour of $\,K_0\,$ we find $\,\|\psi\|^2 \ge
\pi\kappa^{-2}(\ln\kappa a)^{-2} (1\!-\!\eps_3)\,$ for a fixed
$\,\eps_3>0\,$ and $\,a\,$ small enough. Hence
   \begin{equation} \label{upper_3}
{L(\psi)\over \|\psi\|^2}\,<\,-\, {\kappa^2\ln\kappa a\over
\pi(1-\eps_3)}\, (Da^3\ln\kappa a +E)\,,
   \end{equation}
where $\,E:= 2\pi(1+\eps_2)\,$ and
$$
D\,:=\, {2\chi'_1(0)^2 C^2\over (2+\eps_1)\sqrt{\mu_1(1)}}\,.
$$
Minimizing the \rhs of (\ref{upper_3}) with respect to $\,\kappa\,$,
we conclude that to fixed positive $\,\eps_1,\,\eps_2\,$ and
$\,\eps_3\in (0,1)\,$ there is a function $\,g\,$ such that
   \begin{equation} \label{upper bound_3}
{L(\psi)\over \|\psi\|^2}\,<\,g(a) \qquad {\rm and} \qquad
g(a)\,\approx\, -\, {1+\eps_2\over 1-\eps_3}\; {1\over a^2}\:
e^{-2E/Da^3}
   \end{equation}
as $\,a\to 0\,$. The upper bound in (\ref{3-1}) follows readily from
(\ref{upper bound_3}); the extension to the nonsymmetric case is
obtained as in \cite{EV1}.

   \begin{remark}
{\rm In fact, one could suppose $\,M=B_1\,$ because the eigenvalue is
pushed up if we reduce the window to a circle contained in $\,M\,$,
and the obtained bound is all the same not optimal as in
Remark~\ref{optimality}. In the rest of the proof we {\em embedd}
$\,M\,$ into a circle leaving the question about relations between
the constants and the geometry of $\,M\,$ to more sophisticated
methods.}
   \end{remark}

{\em The lower bound} can again be proven in the symmetric case only.
We begin with auxiliary results. When constructing the trial
function component (\ref{trial F_3}), we have used implicitly the
fact that the functional $\,F\,:\; F(\phi)= \int_a^{\infty} \left(
\phi'(t)^2\!+ m^2\phi(t)^2 \right)\, tdt\,$ on $\,C^2([a,\infty))\,$
with the condition $\,\phi(a)=\alpha\,$ and fixed positive $\,a,m\,$
and is minimized by
   \begin{equation} \label{Macdonald}
\phi_0\,:\; \phi_0(t)\,=\, \alpha\, {K_0(mt)\over K_0(ma)}\,,
   \end{equation}
as can be easily seen from solution of the appropriate Euler's
equation. Furthermore, a two--dimensional analogy of the bound (4.2)
in \cite{EV1} is given by the {\em Friedrichs inequality}
\cite[Thm.~1.9]{Ne}: if $\,\Omega\subset\R^n,\; n\ge 2\,$, is a
bounded domain with Lipschitz boundary, there is a positive $\,c\,$
such that
   \begin{equation} \label{Lipschitz}
\|\nabla f\|^2 \ge\, c\|f\|^2
   \end{equation}
holds for every $\,f\in H_0^1(\Omega)\,$. The constant is, of course,
easy to find for the circle $\,\Omega=B_a\,$ in terms of the
appropriate Bessel zero, $\,c=j_{0,1}^2a^{-2}$.

Repeating the argument of \cite{EV1} and the previous section, we
infer that one has to find a lower bound to $\,L(\psi)/\|\psi\|^2$
over all real $\,\psi\in L^2(\Sigma)\,$, which are $\,C^2$, radially
symmetric, and vanish at the boundary except in the window. We can
express such a $\,\psi\,$ in the form of the series (\ref{Fourier})
again, where the convergence is uniform for $\,|x|\ge a\,$. The
coefficients $\,c_n\,$ depend in fact only of $\,r:=|x|\,$. Moreover,
in analogy with (\ref{subexponential decay}) we may restrict our
attention to trial functions with
   \begin{equation} \label{dominated decay_3}
|c_n(r)|\,\le\,|c_n(a)|\, {K_0\left( {\pi\over d}\, \sqrt{n^2\!-1}\,
r \right) \over K_0\left( {\pi\over d}\, \sqrt{n^2\!-1}\, a \right)}
   \end{equation}
for $\,n\ge 2\,$. As before we introduce
   \begin{equation} \label{F_3}
F(x,y)\,:=\, \left\lbrace\; \begin{array}{lll} \alpha\chi_1(y) \quad &
\dots & \quad 0\le r\le 2a \\ c_1(r)\chi_1(y) \quad & \dots & \quad
r\ge 2a \end{array} \right.
   \end{equation}
with $\,\alpha:=c_1(2a)\,$, and divide the rest $\,G(x,y)=
\psi(x,y)-F(x,y)\,$ into
$$
G_1(x,y)\,:=\, (c_1(r)-\alpha)\chi_1(y)
$$
supported in the extended window region, $\,r\le 2a\,$, and
$\,G_2(x,y)= \hat G(x,y)\!+\!\Gamma(x,y)\,$ with
$$
\Gamma(x,y)\,:=\, \sum_{n=2}^{\infty} c_n(2a)\chi_n(y)\,.
$$
We start estimating the reduced energy functional
   \begin{eqnarray} \label{L_3}
L(\psi) \,=\, \|\nabla_x\psi\|^2 +\|G_y\|^2- \left(\pi\over
d\right)^2\|G\|^2- 2\alpha\chi'_1(0) \int_{B_a} G(x,0)\,dx
   \end{eqnarray}
from the ``external" contribution to the first ``two and a half"
terms,
   \begin{eqnarray} \label{L_1}
L_1\!&:=&\! {1\over 2}\,\|\nabla_x\psi\|^2_{r\ge a} +\|G_y\|^2_{r\ge
a}- \left(\pi\over d\right)^2\|G\|^2_{r\ge a} \nonumber \\ \nonumber \\
\!&=&\! \pi\, \sum_{n=1}^{\infty}\, \int_a^{\infty} \left( c'_n(r)^2+\,
2\left(\pi\over d\right)^2(n^2\!-1) c_n(r)^2 \right)\, r\,dr
\nonumber \\ \\
\!&\ge&\! \pi\, \sum_{n=2}^{\infty}\, \int_a^{\infty} \left(
c'_n(r)^2+\, 2\left(\pi n\over d\right)^2 c_n(r)^2 \right)\, r\,dr
\nonumber \\ \nonumber \\
\!&\ge&\! \pi\, \sum_{n=2}^{\infty} c_n(a)^2\, {\pi n\over d}\, a\,
{K_1\left( {\pi n\over d}\, a \right) \over K_0\left( {\pi n\over
d}\, a \right)} \,\ge\, {\pi^2 a\over d}\, \sum_{n=2}^{\infty}
n\,c_n(a)^2\,, \nonumber
   \end{eqnarray}
where in the last line we have used (\ref{Macdonald}), evaluated the
integral as in the first part of the proof, and employed the
inequality $\,K_1(\xi)\ge\,K_0(\xi)\,$ which follows from the
well--known integral representation \cite[9.6.24]{AS}. Next we turn
to
   \begin{equation} \label{L_2}
L_2\,:= \,\|\nabla_x\psi\|^2_{r\le 2a} \,= \, \|\nabla_x
G_1\|^2_{r\le 2a}+ \|\nabla_x G_2\|^2_{r\le 2a}\,.
   \end{equation}
By assumption, $\,G_1\,$ vanishes at $\,r=2a\,$, so the first term
can be estimated from (\ref{Lipschitz}) as
   \begin{equation} \label{G_1}
\|\nabla_x G_1\|^2_{r\le 2a}\,\ge\,{C_1\over 4a^2}\,
\|G_1\|^2_{r\le 2a}\,=\, {C_1\over a^2}\, \|G_1\|^2\,,
   \end{equation}
where $\,4C_1:=j_{0,1}^2\,$. Furthermore, introducing the window
neighbourhood $\,\Omega_a:= B_{2a}\times [0,a]\,$, we have
   \begin{eqnarray} \label{G_2}
\|\nabla_x G_2\|^2_{r\le 2a}\!&=&\!\|\nabla_x \hat G\|^2_{r\le
2a}\,\ge\, {C_1\over a^2}\,\|\hat G\|^2_{r\le 2a} \nonumber \\
\\
\ge\,{C_1\over a^2}\,\|\hat G\|^2_{\Omega_{2a}} \!&\ge&\! {\delta
C_1\over a^2}\,\|\hat G\|^2_{\Omega_{2a}} \,\ge\, {\delta
C_1\over 2a^2}\,\|G_2\|^2_{\Omega_{2a}}-\, {\delta
C_1\over a^2}\,\|\Gamma\|^2_{\Omega_{2a}} \nonumber
   \end{eqnarray}
for all $\,a\le d\,$ and $\,\delta\in (0,1]\,$. The last norm can be
estimated as in the previous cases by combining the smallness of the
$\,\chi_n\,$ norm restricted to $\,[0,a]\,$ with the dominated decay
(\ref{dominated decay_3}),
   \begin{eqnarray*}
\|\Gamma\|^2_{\Omega_a} \!&=&\! 4\pi a^2 \int_0^a \left(
\sum_{n=2}^{\infty} c_n(2a)\chi_n(y) \right)^2 dy \\ \\
&\!\le\!& 8\pi a^2 \left( \sum_{n=2}^{[a^{-1}]+1} n^{-1}
c_n(a)^2 \int_0^a \chi_n(y)^2 dy \right) \sum_{n=2}^{[a^{-1}]+1} n
\\ \\
&\!+\!& 8\pi a^2 \left( \sum^{\infty}_{2\le n=[a^{-1}]+2} n
c_n(a)^2 \int_0^a \chi_n(y)^2 dy \right)
\sum^{\infty}_{2\le n=[a^{-1}]+2} {K_0^2\left( {2\pi a\over d}\,
\sqrt{n^2\!-1} \right) \over n K_0^2\left( {\pi a\over d}\,
\sqrt{n^2\!-1} \right)} \\ \\
&\!\le\!& {16\pi a^3\over d}\, \left( {2\pi^2\over 3d^2}\,+\,
\sum^{\infty}_{2\le n=[a^{-1}]+2} {K_0^2\left( {2\pi a\over d}\,
\sqrt{n^2\!-1} \right) \over n K_0^2\left( {\pi a\over d}\,
\sqrt{n^2\!-1} \right)} \right)\, \sum^{\infty}_{n=2} n c_n(a)^2\,.
   \end{eqnarray*}
The sum in the bracket can be estimated as
$$
\sum^{\infty}_{2\le n=[a^{-1}]+2} {K_0^2\left( {2\pi a\over d}\,
\sqrt{n^2\!-1}\, \right) \over n K_0^2\left( {\pi a\over d}\,
\sqrt{n^2\!-1}\, \right)}\,\le\,
\int_{a^{-1}}^{\infty} {K_0^2\left( {2\pi a\over d}\,
\sqrt{\xi^2\!-1}\, \right) \over \xi K_0^2\left( {\pi a\over d}\,
\sqrt{\xi^2\!-1}\, \right)} \,d\xi \,\le\,
\int_1^{\infty} {K_0^2\left( {\pi\xi\over d} \right) \over \xi
K_0^2\left( {\pi\xi\over 2d} \right)}
$$
for $\,a< {1\over 2}\sqrt{3}\,$, and the integral on the \rhs is
convergent, because $\,K_0(\xi)\approx \sqrt{\pi\over 2\xi}\,
e^{-\xi}\,$ as $\,\xi\to\infty\,$. Hence there is a positive
$\,C_2\,$ independent of $\,\psi\,$ and $\,a\,$ such that
   \begin{equation} \label{Gamma_3 norm}
{C_1\over a^2}\:\|\Gamma\|^2_{\Omega_a}\,<\, C_2 a
\sum^{\infty}_{n=2} n c_n(a)^2\,.
   \end{equation}
Combining the estimates (\ref{L_1})--(\ref{Gamma_3 norm}), we arrive
at
$$
L_1\!+L_2 \,\ge\, a\left({\pi^2\over d}\,-\delta C_2\right)
\sum^{\infty}_{n=2} n c_n(a)^2+\, {C_1\over a^2}\, \|G_1\|^2 +\,
{\delta C_1\over 2a^2}\, \|G_2\|^2_{\Omega_a}\,,
$$
which gives
   \begin{equation} \label{L_12}
L_1\!+L_2 \,\ge\, {C_1\over a^2}\, \|G_1\|^2 +\,{m^2\over a^2}\,
\|G_2\|^2_{\Omega_a}
   \end{equation}
for some $\,m>0\,$ and all sufficiently small $\,a\,$.

The norm of $\,G_y\,$ is estimated as in the two--dimensional case
\cite{EV1},
$$
\|G_y\|^2_{r\le a}\,\ge\, \|G_{2,y}\|^2_{r\le a} -\,{2\pi\over d^2}\,
\left( 2\|G_1\|^2_{r\le a}\!+ d\|G_2(\cdot,0)\|^2_{r\le a}\right)\,,
$$
which together with (\ref{L_12}) yields
   \begin{eqnarray*}
L_1 \!\!&+&\!L_2\,+\:\|G_y\|^2_{r\le a}-\, \left(\pi\over d\right)^2
\|G\|^2_{r\le a} \\ \\
\!&\ge &\! \|G_{2,y}\|^2_{r\le a}-\, \left(\pi\over d\right)^2
\|G_2\|^2_{r\le a}+\, \left( {C_1\over a^2}\,-\,{\pi(\pi+4)\over d^2}
\right) \|G_1\|^2_{r\le a} \\ \\
\!&+&\! {m^2\over a^2}\, \|G_2\|^2_{\Omega_a}-\, {2\pi\over d}\,
\|G_2(\cdot,0)\|^2_{r\le a} \\ \\
\!&\ge &\! \left( {c_0\over a}\,-\,{2\pi\over d}\right)
\|G_2(\cdot,0)\|^2_{r\le a} \,\ge\, {c_0\over
2a}\,\|G_2(\cdot,0)\|^2_{r\le a}
   \end{eqnarray*}
for a positive $\,c_0\,$ and any $\,a\,$ small enough; in the second
step we have neglected a positive term and employed
Lemma~\ref{lemma 3}. Substituting from here to (\ref{L_3}) and using
the Schwarz inequality,
$$
\int_{B_a} G(x,0)\,dx \,\le\, \|G_2(\cdot,0)\|_{r\le a}
\sqrt{\pi}\,a\,
$$
we get
   \begin{eqnarray*}
L(\psi) \!&\ge&\! {1\over 2}\, \|\nabla_x\psi\|^2_{r\ge 2a}
-\,2\alpha a\chi'_1(0) \sqrt{\pi}\, \|G_2(\cdot,0)\|_{r\le a}+\,
{c_0\over 2a}\, \|G_2(\cdot,0)\|_{r\le a}^2 \\ \\
\!&\ge&\! {1\over 2}\, \|\nabla_x\psi\|^2_{r\ge 2a} -\,{2\pi\alpha^2
\chi'_1(0)^2\over c_0}\, a^3 \,.
   \end{eqnarray*}
The first term on the \rhs can be estimated from below by the first
transverse--mode contribution. The same applies to $\,\|\psi\|^2$, so
finally we find
   \begin{equation} \label{lbound_3}
{L(\psi)\over \|\psi\|^2}\,\ge\, {\int_{2a}^{\infty} c'_1(r)^2 r\,dr
-\,{\pi\chi'_1(0)^2\over c_0}\, a^3 c_1(2a)^2 \over 2
\int_{2a}^{\infty} c_1(r)^2 r\,dr}\,.
   \end{equation}
In analogy with (\ref{Macdonald}) one has to solve the appropriate
Euler's equation to check that the \rhs of (\ref{lbound_3}) is
minimized by $\,c_1=\phi_{\kappa}$ for some $\,\kappa>0\,$, where
$$
\phi_{\kappa}(r)\,:=\, c_1(2a)\, {K_0(\kappa r)\over
K_0(2\kappa a)}\,.
$$
Substituting into (\ref{lbound_3}), evaluating the integrals, and
taking the asymptotics for small $\,a\,$, we infer that
$$
{L(\psi)\over \|\psi\|^2}\,\ge\, -\,\kappa^2 \ln(2\kappa a)\, \left(
{\pi\chi'_1(0)^2 \over c_0(1+\eps_2)}\, a^3\ln(2\kappa a)+
\,{1-\eps_1\over 1+\eps_2} \right)
$$
holds for any fixed $\,\eps_1,\,\eps_2>0\,$ and all sufficiently
small $\,a\,$. It remains to find the minimum of the \rhs with
respect to $\,\kappa\,$. However, since it differs from
(\ref{upper_3}) just by the values of the constants, the argument
is concluded as in the first part of the proof. \quad \QED

\section{Conclusions}
\setcounter{equation}{0}

To make sense of the derived bounds one has to take into account two
aspects of the problem. First of all, we have mentioned already that
the discrete spectrum can also be found numerically by means
of the mode--matching method; a detailed description of the
two--dimensional case is given in \cite{ESTV}. Although the method
converges rather slowly if the window is narrow, the results obtained
for a single window clearly suggest that the true asymptotics exists
and is of the same type as our asymptotic bounds.

Another insight can be obtained from comparing our result with the
well-known weak--coupling asymptotics for Schr\"odinger operators in
dimension one and two \cite{BGS,Kl,Si}. The ground state of the
coupled strips in the narrow--window case is dominated the lowest
transverse--mode component with long exponentially decaying tails and
a local modification in the coupling region. In a similar way, a link
can be made between window--connected layers and a two--dimensional
Schr\"odinger operator. The comparison shows that the attractive
interaction due to opening a narrow window (in particular, by
changing the Dirichlet b.c. to Neumann at a short segment of the
boundary in the symmetric case) acts effectively as a potential well
of a depth proportional to the size of the window.

   \begin{conjecture}
Let $\,H(d_1,d_2;\WW)\,$ be the operators described above. The
ground--state eigenvalue behaves for small $\,|\WW|\,$ as
   \begin{eqnarray}
\epsilon(a) \,\approx\, \left(\pi\over d\right)^2 -\,{1\over d^2}\,
\left(\sum_{k=1}^N c_{2,k}(\nu) a_k^2 \right)^2 & \quad \dots \; &
\dim\Sigma=2 \\ \nonumber \\
\epsilon(a) \,\approx\, \left(\pi\over d\right)^2 -\,{1\over
d^2}\,\exp\left\lbrace -\left(\sum_{k=1}^N c_{3,k}(\nu) a_k^3
\right)^{-1}\right\rbrace  & \quad \dots \; & \dim\Sigma=3
   \end{eqnarray}
where $\,\nu:=d^{-1}\min\{d_1,d_2\}\,$, and $\,a_k\,$ in the
three-dimensional case is the scaling parameter of the $\,k$--th
window.
   \end{conjecture}
The conjecture is based on the described analogy only, and therefore
it is difficult to say more about the coefficients. It is not
excluded that they depend on the geometry of the window--center set
for $\,N>1\,$; in the three-dimensional case the shapes of the scaled
windows may also play role. We refrain from speculating about
the nature of the error terms.

On the other hand, we are convinced that the open ``constant
cross--section" shape of our regions $\,\Sigma\,$ is crucial for the
asymptotics. For instance, if $\,\Sigma\,$ is instead a bounded
planar region with the Dirichlet boundary in which we open a window
(to another bounded region the essential spectrum threshold of which
is not lower) or a Neumann segment, we conjecture that leading term
in the ground state shift is proportional to the {\em square} of the
window width. Moreover, the same asymptotics is expected to be valid
for higher eigenvalues provided the corresponding eigenfunctions are
locally {\em symmetric} with respect to the window axis.  In any
case, proving of such asymptotic properties represents an intriguing
mathematical problem.

\section*{Acknowledgments}

Thanks are due to P.~Lindovsk\'{y} for useful discussions. The
research has been partially supported by the Grant GACR No.
202--0218.

\end{document}